\documentclass[11pt,twoside,reqno]{amsart}
\usepackage{amsmath,amsfonts,amsthm,amssymb}
\usepackage{graphicx,psfrag,overpic}
\usepackage{bm}
\usepackage[all,cmtip,line]{xy}
\usepackage{array}
\newtheorem{problemL}{Problem}

\theoremstyle{definition}

\newcommand \gam{\gamma}
\newcommand \R{\mathbb{R}}

\newcommand \der{\partial}
\newcommand \vphi{\varphi}

\newcommand \alp{\alpha}

\newcommand \irho{\rho_{\infty}}

\newcommand{\PtIncW}{{P_0}}
\newcommand{\PtUpL}{{P_1}}
\newcommand{\PtLwL}{{P_2}}
\newcommand{\PtLwR}{{P_3}}
\newcommand{\PtUpR}{{P_4}}

\newlength{\originalbase}
\begin{document}

\title[Shock Reflection-Diffraction and Nonlinear Equations of Mixed Type]
{Shock Reflection-Diffraction,\\ von Neumann's Conjectures, \\and Nonlinear Equations of Mixed Type}

\author{Gui-Qiang Chen}
\address{Gui-Qiang G. Chen, Mathematical Institute, University of Oxford,
Oxford, OX2 6GG, England;
School of Mathematical Sciences, Fudan University, Shanghai  200433, China}
\email{chengq@maths.ox.ac.uk}

\author{Mikhail Feldman}
\address{Mikhail Feldman, Department of Mathematics\\
         University of Wisconsin\\
         Madison, WI 53706-1388, USA}
\email{feldman@math.wisc.edu}
\maketitle

\begin{abstract}
Shock waves are fundamental in nature.
One of the most fundamental problems in fluid mechanics is shock reflection-diffraction by wedges.
The complexity of reflection-diffraction configurations was first reported
by Ernst Mach in 1878.
The problems remained dormant until the 1940s when John von Neumann,
as well as other mathematical/experimental
scientists, began extensive research into all aspects of shock reflection-diffraction phenomena.
In this paper we start with shock reflection-diffraction phenomena
and historic perspectives, their fundamental scientific issues
and theoretical roles
in the mathematical theory of hyperbolic systems of conservation laws.
Then we present how the global shock reflection-diffraction problem
can be formulated as a boundary value problem in an unbounded domain
for nonlinear conservation laws of mixed hyperbolic-elliptic type,
and describe the von Neumann conjectures: the sonic conjecture and
the detachment conjecture.
Finally we discuss some recent developments in solving the von Neumann conjectures
and establishing a mathematical theory of shock reflection-diffraction,
including the existence, regularity, and stability of global regular
configurations of shock reflection-diffraction by wedges.
\end{abstract}

\bigskip

\section{Introduction}

\numberwithin{equation}{section}

Shock waves are steep fronts that propagate in the compressible fluids when
the convective motion dominates the diffusion. They are fundamental in
nature. Examples include shock waves formed by solar winds (bow shocks),
supersonic or near sonic aircrafts (transonic shocks around the body or
supersonic bubbles with transonic shocks), explosions (blast waves), and
various other natural processes. When a shock hits an obstacle (steady or
flying), shock reflection-diffraction phenomena occur.
One of the most fundamental
problems in fluid mechanics is the problem of shock
reflection-diffraction by wedges; see Ben-Dor \cite{BD},
Courant-Friedrichs \cite{CF}, von Neumann \cite{Neumann1,Neumann2,Neumann},
and the references cited therein. When a plane
shock hits a wedge head on, it experiences a reflection-diffraction process,
and then a fundamental question is what types of wave patterns of
reflection-diffraction configurations may be formed around the wedge.

\smallskip
The complexity of reflection-diffraction configurations was first reported
by Ernst Mach \cite{Mach} in 1878, who first observed two patterns of
reflection-diffraction configurations:
regular reflection (two-shock configuration; see Fig. 1 (left))
and Mach reflection (three-shock/one-vortex-sheet configuration;
see Fig. 1 (center)); also see \cite{BD,Chen-F-7, CF, VD}. The issues
remained dormant until the 1940¡¯s when John von Neumann,
as well as other mathematical/experimental scientists,
began extensive research into all aspects of shock reflection-diffraction phenomena,
due to its importance in applications.
See von Neumann \cite{Neumann1,Neumann2,Neumann},
Courant-Friedrichs \cite{CF}, Glimm-Majda \cite{GlimmMajda},
and Ben-Dor \cite{BD};
also see \cite{Chen-F-7,VD} and the references cited therein.
It has been found that the situations are much more complicated than
what Mach originally observed: The Mach reflection can be further
divided into more specific sub-patterns, and various other patterns of
shock reflection-diffraction may occur such as the double Mach reflection,
von Neumann reflection, and Guderley reflection;
see \cite{BD,Chen-F-7,CF,GlimmMajda,VD}
and the references cited therein.
The fundamental scientific issues include the following:

\smallskip
\begin{itemize}
\item[(i)] Structure of the shock reflection-diffraction configurations;

\item[(ii)] Transition criteria between the different patterns of shock
reflection-diffraction configurations;

\item[(iii)] Dependence of the patterns upon the physical parameters such as the
wedge angle $\theta_w$, the incident-shock-wave Mach number, and the
adiabatic exponent $\gamma\ge 1$.
\end{itemize}

\smallskip
In particular, several transition criteria between the different
patterns of shock reflection-diffraction configurations have been proposed,
including the sonic conjecture and the detachment conjecture by von Neumann
\cite{Neumann1,Neumann2,Neumann}.

\begin{figure}[h]
\centering
\includegraphics[height=1.3in,width=6.0in]{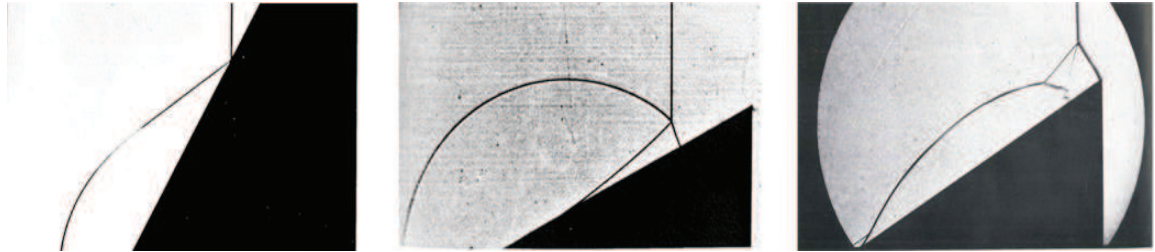} 
\caption[]{Regular reflection (left); $\,\,\,$ Mach reflection (center); $\,\,\,$ Irregular Mach reflection (right).
From  Van Dyke \cite{VD}, pages 142--144}
 \label{fig:Experiment}
\end{figure}

Careful asymptotic analysis has been made for various reflection-diffraction
configurations in Lighthill \cite{Lighthill1,Lighthill2},
Keller-Blank \cite{KB}, Hunter-Keller \cite{HK}, Harabetian \cite{Harabetian},
Morawetz \cite{Morawetz2}, and the references
cited therein; also see Glimm-Majda \cite{GlimmMajda}.
Large or small scale numerical simulations have been also made;
{\it cf.} \cite{BD,GlimmMajda,WC} and the references cited therein.
However, most of the fundamental issues for shock reflection-diffraction
phenomena have not been understood, especially the global structure and
transition between different patterns of shock reflection-diffraction configurations.
This is partially because physical and numerical experiments are
hampered by various difficulties and have not yielded clear transition criteria
between different patterns. In particular, numerical dissipation or physical
viscosity smear the shocks and cause boundary layers that interact with
the reflection-diffraction patterns and can cause spurious Mach steams; {\it cf.}
Woodward-Colella \cite{WC}. Furthermore, some different patterns occur when
the wedge angles are only fractions of a degree apart, a resolution even by
sophisticated modern experiments ({\it cf.} \cite{LD}) has been unable to reach.
For this reason, it is almost impossible to
distinguish experimentally between the sonic and detachment criteria,
as pointed out in \cite{BD}.
In this regard, the necessary approach to understand
fully the shock reflection-diffraction phenomena, especially the transition criteria,
is still via rigorous mathematical analysis.
To achieve this, it is essential to establish first the global existence, regularity,
and structural stability of solutions of the shock reflection-diffraction problem.

\smallskip
Furthermore, shock reflection-diffraction configurations are the core configurations
in the structure of global entropy solutions of the two-dimensional
Riemann problem for hyperbolic conservation laws,
while the Riemann solutions are building blocks and local structure of general
solutions and determine global attractors and asymptotic states of entropy
solutions, as time tends to infinity, for multidimensional hyperbolic systems
of conservation laws. See \cite{Chen-F-7,GlimmMajda,Serre,Zhe}
and the references cited therein. In this sense, we have to understand
the shock reflection-diffraction phenomena, in order to understand fully
global entropy solutions to multidimensional hyperbolic systems of conservation
laws.

\smallskip
In Section 2, we first formulate the shock reflection-diffraction problem
into an initial-boundary value problem.
Based on the invariance of the problem under the self-similar scaling
and some basic properties
of the governing equations,
we reformulate the problem as a boundary value problem in an unbounded domain
for nonlinear conservation laws of mixed hyperbolic-elliptic type.
In Section 3, we describe the von Neumann's conjectures: the sonic conjecture
and the detachment conjecture.
In Section 4, we
discuss some recent developments in solving the von Neumann's conjectures
and establishing a mathematical theory of shock reflection-diffraction,
including the existence, regularity, and stability of global regular configurations
of shock reflection-diffraction by wedges.
The shock reflection-diffraction problems involve several core mathematical difficulties
we have to face in solving nonlinear partial differential equations in mechanics and geometry;
also see \cite{CSW1,CSW2,HH,Lin,Yau} and the references cited therein.
These include nonlinear equations of mixed hyperbolic-elliptic type,
nonlinear degenerate elliptic equations, nonlinear degenerate hyperbolic equations,
free boundary problems for nonlinear degenerate equations, corner singularity/regularity especially
when free boundaries meet degenerate curves, and {\it a priori} estimate techniques.

\section{Mathematical Formulation and Nonlinear Equations of Mixed Type}
\label{section:shock-polar}

\numberwithin{equation}{section}

In this section we first formulate the shock reflection-diffraction problem
into an initial-boundary value problem and then
reformulate the problem as a boundary value problem in an unbounded domain for a
nonlinear conservation laws of mixed elliptic-hyperbolic type.

\subsection{Mathematical Problems}
The potential flow is governed by the conservation law of mass and the Bernoulli law:
\begin{align}
\label{1-a}
&\der_t\rho+ \nabla_{\bf x}\cdot (\rho \nabla_{\bf x}\Phi)=0,\\
\label{1-b}
&\der_t\Phi+\frac 12|\nabla_{\bf x}\Phi|^2+h(\rho)=B
\end{align}
for the density $\rho$ and the velocity potential $\Phi$,
where $B$ is the Bernoulli constant determined by the incoming flow and/or boundary conditions, and
\begin{equation*}
h'(\rho)=\frac{p'(\rho)}{\rho}=\frac{c^2(\rho)}{\rho}
\end{equation*}
with $c(\rho)$ being the sound speed and $p$ the pressure.
For an ideal polytropic gas, the pressure $p$ and the sound speed $c$ are given by
\begin{equation*}
p(\rho)=\kappa \rho^{\gam}, \quad c^2(\rho)=\kappa \gam \rho^{\gam-1}
\end{equation*}
for constants $\gam>1$ and $\kappa>0$.
By scaling, we may choose $\kappa=1/\gam$ to have
\begin{equation}
\label{1-c}
h(\rho)=\frac{\rho^{\gam-1}-1}{\gam-1},\quad c^2(\rho)=\rho^{\gam-1}.
\end{equation}

\smallskip
From \eqref{1-b} and \eqref{1-c}, we have
\begin{equation}
\label{1-b1}
\rho(\der_t\Phi,\nabla_{\bf x}\Phi)=h^{-1}\bigl(B-\der_t\Phi-\frac 12|\nabla_{\bf x}\Phi|^2\bigr).
\end{equation}
Then system \eqref{1-a}--\eqref{1-b} can be rewritten as
\begin{equation}
\label{1-b2}
\der_t\rho(\der_t\Phi, \nabla_{\bf x}\Phi)
+\nabla_{\bf x}\cdot\big(\rho(\der_t\Phi, \nabla_{\bf x}\Phi)\nabla_{\bf x}\Phi\big)=0
\end{equation}
with $\rho(\der_t\Phi, \nabla_{\bf x}\Phi)$ determined by \eqref{1-b1}.

\smallskip
When a vertical planar shock perpendicular to the flow
direction $x_1$ and separating two uniform states (0) and (1), with constant velocities
$(u_0, v_0)= (0, 0)$ and $(u_1,v_1)=(u_1, 0)$,
and constant densities $\rho_1>\rho_0$
(state (0) is ahead or to the right of shock in the figure and state
(1) is behind the shock), hits a symmetric wedge:
$$
W:= \{|x_2| < x_1 \tan(\theta_w), x_1 > 0\}
$$
head on at time $t = 0$, then
a reflection-diffraction process takes place when $t > 0$.
Since state (1) does not satisfy the slip
boundary condition $\partial_\nu\varphi= 0$  prescribed on the wedge boundary,
the solution must differ from
state (1) on and near
the wedge boundary.
Mathematically, the shock reflection-diffraction problem
is a multidimensional
lateral
Riemann problem in the domain
$\R^2\setminus \bar{W}$.

\medskip
\begin{problemL}[Lateral Riemann Problem]\label{problem-1}
{\it
Piecewise constant initial data, consisting of state $(0)$
on $\{x_1>0\}\setminus \bar{W}$
and state $(1)$ on $\{x_1 < 0\}$, are prescribed at $t = 0$.
Seek a solution of $(1)$ for $t\ge 0$ subject to these initial
data and the boundary condition $\nabla\Phi\cdot\nu=0$ on $\partial W$.
}
\end{problemL}

\medskip
Notice that the initial-boundary value problem for \eqref{1-a}--\eqref{1-b}
is invariant under the scaling:
\begin{equation*}
(t, {\bf x})\rightarrow (\alp t, \alp{\bf x}),\quad (\rho, \Phi)\rightarrow (\rho, \frac{\Phi}{\alp})
\qquad \text{for}\;\;\alp\neq 0.
\end{equation*}
Thus, we seek self-similar solutions in the form of
\begin{equation*}
\rho(t, {\bf x})=\rho(\xi,\eta),\quad \Phi(t, {\bf x})=t\phi(\xi,\eta)\qquad\text{for}\;\;(\xi,\eta)=\frac{{\bf x}}{t}.
\end{equation*}
Then the pseudo-potential function $\vphi=\phi-\frac 12(\xi^2+\eta^2)$ satisfies the following
nonlinear conservation laws of mixed type:
\begin{equation}
\label{2-1}
{\rm div}\big(\rho(|D\vphi|^2,\vphi)D\vphi\big)+2\rho(|D\vphi|^2,\vphi)=0,
\end{equation}
for
\begin{equation}
\label{2-1a}
\rho(|D\vphi|^2,\vphi)=
\bigl(\rho_0^{\gamma-1}-(\gam-1)(\frac 12|D\vphi|^2+\vphi)\bigr)^{\frac{1}{\gam-1}},
\end{equation}
where the divergence ${\rm div}$ and gradient $D$ are with respect to $(\xi,\eta)$,
and we have used the initial condition to fix the Bernoulli constant $B$ in \eqref{1-b}.
Then equation \eqref{2-1} is a second-order nonlinear conservation law of mixed elliptic-hyperbolic type.
It is elliptic if and only if
\begin{equation}
\label{1-f}
|D\vphi|^2<c^2(|D\vphi|^2,\vphi):=\rho_0^{\gamma-1}-(\gam-1)(\frac 12|D\vphi|^2+\vphi),
\end{equation}
which is equivalent to
\begin{equation}
\label{1-g}
|D\vphi|<c_*(\vphi,\irho):=\sqrt{\frac{2}{\gam+1}\bigl(\rho_0^{\gamma-1}-(\gam-1)\vphi\bigr)}.
\end{equation}
The type of equation, due to its
nonlinearity, depends on the solution of the equation, which
makes the problem truly fundamental and challenging.
Elliptic regions of the solution correspond to the subsonic flow, and hyperbolic
regions to the supersonic flow.

\smallskip
Shocks are discontinuities in the pseudo-velocity $D\varphi$.
If curve $S$ subdivides the domain into subdomains $\Omega^\pm$, and
$\varphi\in C^1(\Omega^\pm\cup S)\cap C^2(\Omega)$,
then $\varphi$ is a weak solution of \eqref{2-1} in $\Omega$ across $S$ if and only if
$\varphi$ satisfies \eqref{2-1} in $\Omega^\pm$
and the Rankine-Hugoniot conditions on $S$:
\begin{equation}\label{RH}
[\varphi]_S=0, \qquad [\rho(|D\varphi|^2, \varphi)D\varphi\cdot\nu]_S=0,
\end{equation}
where $[\cdot]_S$ is the jump across $S$, and $\nu$ is the unit normal to $S$.

\smallskip
Then Problem \ref{problem-1} is reformulated as the following boundary value problem
in unbounded domain
$$
\Lambda:=\R^2\setminus\{|\eta|\le \xi \tan(\theta_w), \xi>0\}
$$
in the self-similar coordinates $(\xi,\eta)$.

\begin{problemL}[Boundary Value Problem]\label{problem-2}
{\it Seek a
solution $\varphi$ of equation \eqref{2-1} in the self-similar
domain $\Lambda$ with the slip boundary condition
$D\varphi\cdot\nu|_{\partial\Lambda}=0$
on the wedge boundary $\partial\Lambda$
and the asymptotic boundary condition at infinity:
$$
\varphi\to\bar{\varphi}=
\begin{cases} \varphi_0 \qquad\mbox{for}\,\,\,
                         \xi>\xi_0, |\eta|>\xi \tan(\theta_w),\\
              \varphi_1 \qquad \mbox{for}\,\,\,
                          \xi<\xi_0.
\end{cases}
\qquad \mbox{when $\xi^2+\eta^2\to \infty$.}
$$
}
\end{problemL}

By symmetry we can restrict to the upper half-plane $\{\eta>0\}\cap\Lambda$,
with condition $\partial_\nu \varphi=0$ on $\{\eta=0\}\cap\Lambda$.

\section{von Neumann's Conjectures on Shock Reflection-Diffraction Configurations}

\numberwithin{equation}{section}

In this section, we discuss von Neumann's conjectures,
including the sonic conjecture and the detachment conjecture, on
shock reflection-diffraction configurations.

\begin{figure}
\centering
\begin{minipage}{0.435\textwidth}
\centering
\psfrag{Om}{\Large$\Omega$}
\psfrag{P0}{$P_0$}
\psfrag{P1}{$P_1$}
\psfrag{P2}{$P_2$}
\psfrag{P3}{$P_3$}
\psfrag{P4}{$P_4$}
\includegraphics[height=1.7in,width=2.8in]{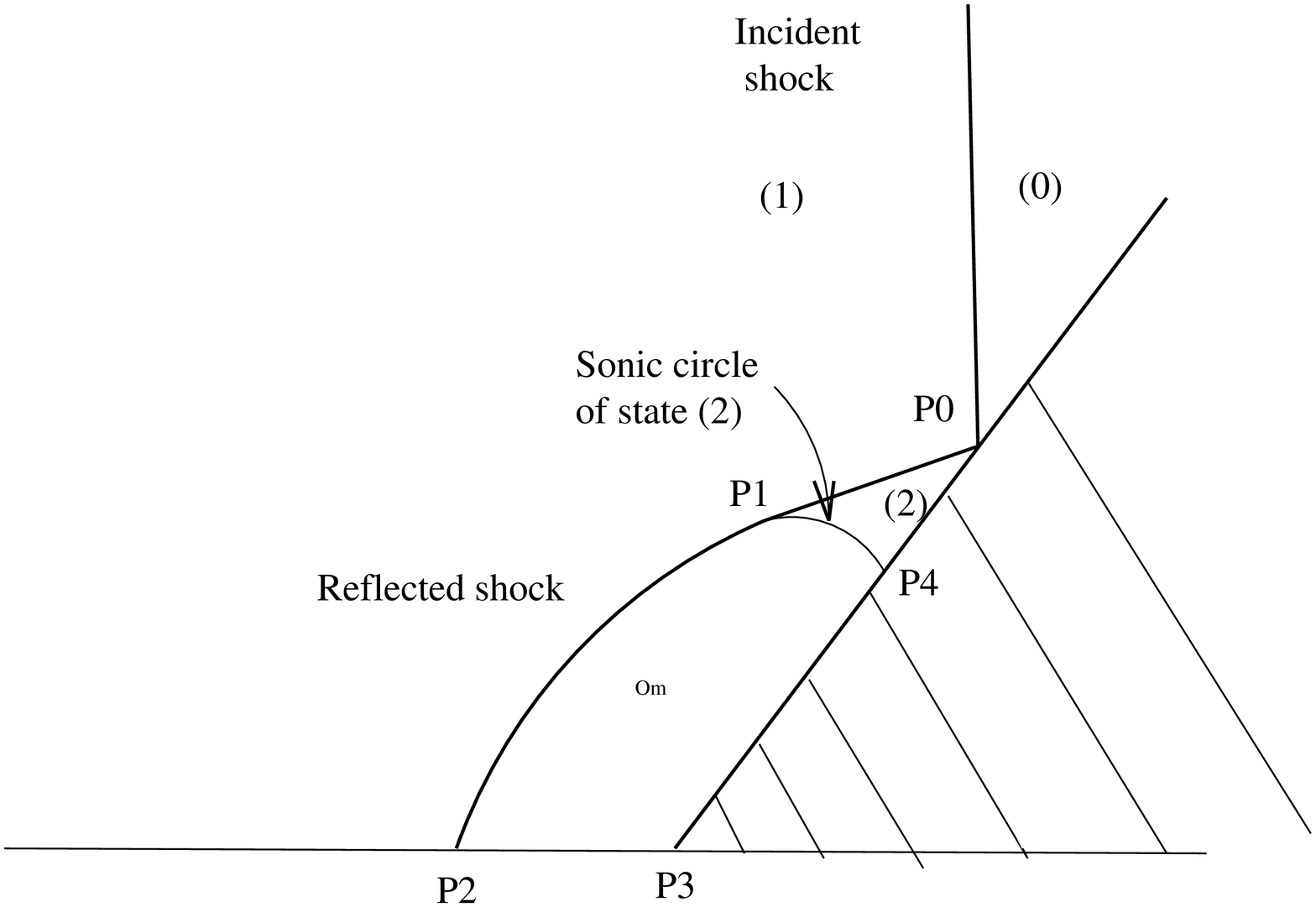}  
\caption{Supersonic regular reflection}
\label{fig:RegularReflection}
\end{minipage}
\qquad
\begin{minipage}{0.435\textwidth}
\centering
\psfrag{Om}{\Large$\Omega$}
\psfrag{p0}{$P_0$}
\psfrag{p1}{$P_1$}
\psfrag{p2}{$P_2$}
\psfrag{p3}{$P_3$}
\psfrag{p4}{$P_4$}
\includegraphics[height=1.7in,width=2.7in]{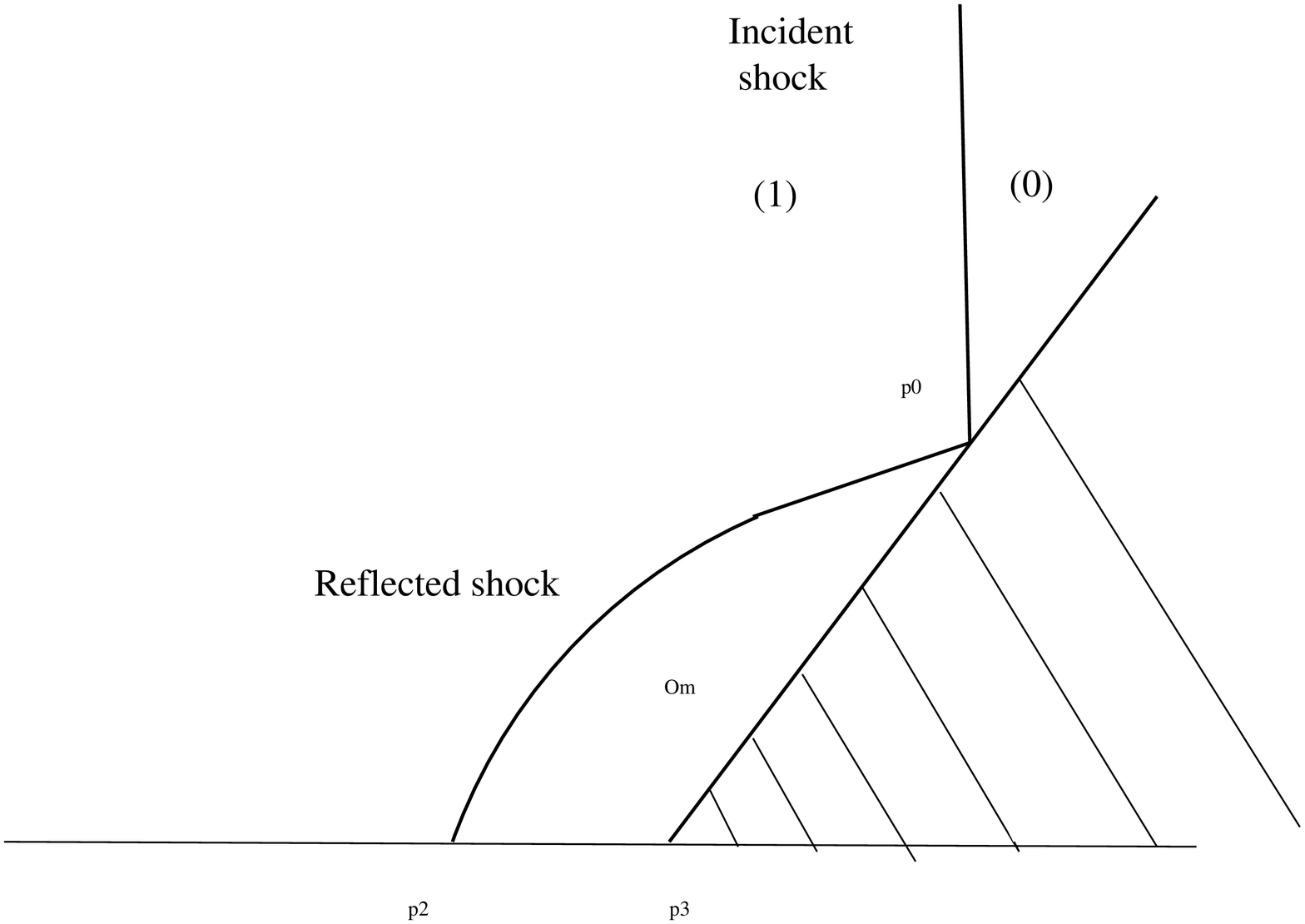}  
\caption{Subsonic $\,\,\,$ regular reflection}\label{fig:SubsoncRegularReflection}
\end{minipage}
\end{figure}

\smallskip
If a solution has the regular reflection-diffraction configurations as shown in Fig. \ref{fig:RegularReflection}
and Fig. \ref{fig:SubsoncRegularReflection},
and if $\varphi$ is smooth in the subregion between the reflected shock and the wedge, then
it should satisfy  the boundary condition $D\varphi\cdot\nu=0$ and the Rankine-Hugoniot
conditions (\ref{RH}) at $\PtIncW$
across the reflected shock separating it from state (1).
We define the uniform state (2) with pseudo-potential $\varphi_2(\xi,\eta)$ such
that
$$
D\varphi_2(\PtIncW)=D\varphi(\PtIncW),
$$
and the constant density $\rho_2$ of state (2) is equal to $\rho(|D\varphi|^2, \varphi)(\PtIncW)$
defined by (\ref{2-1}):
$$
\rho_2=\rho(|D\varphi|^2, \varphi)(\PtIncW).
$$
Then $D\varphi_2\cdot\nu=0$ on the wedge boundary, and the Rankine-Hugoniot conditions
(\ref{RH}) hold on the flat shock $S_1=\{\varphi_1=\varphi_2\}$
between states (1) and (2), which passes through $\PtIncW$.

\smallskip
State (2) can be either subsonic or supersonic at $\PtIncW$.
This determines the subsonic or supersonic type of regular reflection-diffraction configurations.
The supersonic regular reflection-diffraction configuration
as shown in Fig. \ref{fig:RegularReflection}
consists of three uniform states (0), (1), (2), and a non-uniform state in the domain $\Omega$,
where the equation is elliptic.
The reflected shock $\PtIncW\PtUpL\PtLwL$ has a straight part $\PtIncW\PtUpL$.
The elliptic domain $\Omega$ is separated from the hyperbolic region $\PtIncW\PtUpL\PtUpR$
of state (2) by a sonic arc $\PtUpL \PtUpR$.
The subsonic regular reflection-diffraction configuration
as shown in Fig. \ref{fig:SubsoncRegularReflection}
consists of two uniform states (0) and (1), and a non-uniform
state in the domain $\Omega$, where the equation is elliptic and
$D(\varphi_{|\Omega})(\PtIncW)=D\varphi_2(\PtIncW)$.

\smallskip
Thus, a necessary condition for the existence of regular reflection
solution is the existence of the uniform state (2) determined by the conditions described above.
These conditions lead to an algebraic system for the constant velocity $(u_2, v_2)$ and
density $\rho_2$ of state (2),
which has solutions for some but not all of the wedge angles.
Specifically, for fixed densities $\rho_0<\rho_1$ of states (0) and (1),
there exist a sonic angle $\theta_w^s$ and a detachment angle $\theta_w^{d}$
satisfying
$$
0<\theta_w^{d}<\theta_w^s<\frac{\pi}{2}
$$
such that state (2) exists for all $\theta_w\in (\theta_w^{d}, \frac{\pi}{2})$
and does not exist for $\theta_w\in (0, \theta_w^{d})$,
and the weak state (2) is

\smallskip
\begin{itemize}
\item[(i)] supersonic at the reflection point $\PtIncW(\theta_w)$ for
$\theta_w\in (\theta_w^s, \frac{\pi}{2})$;

\item[(ii)] sonic for $\theta_w=\theta_w^s$;

\item[(iii)]
subsonic for $\theta_w\in (\theta_w^{d}, \hat\theta_w^s)$
\end{itemize}
for some $\hat\theta_w^s\in(\theta_w^{d}, \theta_w^s]$.
In fact, for each $\theta_w\in(\theta_w^{d}, \frac{\pi}{2})$,
there exists also a {\em strong} state (2) with $\rho_2^{strong}>\rho_2^{weak}$.
There had been a long debate to
determine which one is physical for the local theory; see
Courant-Friedrichs \cite{CF}, Ben-Dor \cite{BD}, and the references
cited therein.
It is expected that strong reflection-diffraction configurations
are non-physical; indeed, it is shown that
the weak reflection-diffraction configuration tends to the unique normal
reflection, but the strong reflection-diffraction configuration does not,
when the wedge angle $\theta_w$ tends
to $\frac{\pi}{2}$, as proved in Chen-Feldman \cite{Chen-F-6}.

\smallskip
If the weak state (2) is supersonic, the propagation speeds of the solution
are finite and state (2) is
completely determined by the local information: state (1), state
(0), and the location of point $P_0$. That is, any information
from the region of reflection-diffraction,
especially the disturbance at corner $P_3$,
cannot travel towards the reflection point $P_0$.
However, if
it is subsonic, the information can reach $P_0$ and interact with
it, potentially altering a different reflection-diffraction configuration.
This argument motivated the following conjecture von Neumann in \cite{Neumann1,Neumann2}:

\medskip
{\bf The Sonic Conjecture}:
{\em There exists a supersonic reflection-diffraction
configuration when $\theta_w\in
(\theta_w^s, \frac{\pi}{2})$ for $\theta_w^s>\theta_w^d$,
{\it i.e.},
the supersonicity of the weak state {\rm (2)} implies the existence
of a supersonic regular reflection-diffraction
solution, as shown in Fig. {\rm 3}.}

\medskip
Another conjecture is that global regular reflection-diffraction
configuration is possible whenever the local regular reflection at the reflection
point is possible:

\medskip
{\bf The Detachment Conjecture}: {\em There
exists a regular reflection-diffraction configuration for
all wedge angles $\theta_w\in (\theta_w^{d}, \frac{\pi}{2})$, {\it i.e.},
that the existence of state {\rm (2)} implies the existence
of a regular reflection-diffraction
solution, as shown in Fig. {\rm 3}.}

\medskip
It is clear that the supersonic/subsonic regular reflection-diffraction configurations are
not possible without a local two-shock configuration at the
reflection point on the wedge, so this is the weakest possible
criterion for the existence of (supersonic/subsonic) regular
shock reflection-diffraction configurations.

\section{Solution to the von Neumann Conjectures}

We observe that the key obstacle to the existence
of regular shock reflection-diffraction configurations as conjectured by von Neumann \cite{Neumann1,Neumann2}
is an additional possibility that,
for some wedge angle $\theta_w^a\in (\theta_w^d, \frac{\pi}2)$, shock
$\PtIncW\PtLwL$ may attach to the wedge tip $\PtLwR$, as observed
by experimental results ({\it cf.} \cite[Fig. 238]{VD}).
To describe the conditions of such attachments, we note that
$$
\rho_1>\rho_0, \qquad
u_1=(\rho_1-\rho_0)
\sqrt{\frac{2(\rho_1^{\gamma-1}-\rho_0^{\gamma-1})}{\rho_1^2-\rho_0^2}}.
$$
Then, for each $\rho_0$, there exists $\rho_1^{cr}>\rho_0$ such that
\begin{eqnarray*}
&& u_1\le c_1 \quad \mbox{if $\rho_1\in (\rho_0, \rho_1^{cr}]$},\\[2mm]
&& u_1>c_1 \quad \mbox{if $\rho_1\in (\rho_1^{cr}, \infty)$}.
\end{eqnarray*}

\smallskip
If $u_1\le c_1$, we can rule out the solutions with shocks attached to the
wedge tip.

\smallskip
If $u_1> c_1$, this is unclear and could not be even true
as experiments show ({\it e.g.} \cite[Fig. 238]{VD}).

\smallskip
Thus, in \cite{Chen-F-6,Chen-F-7},
we have obtained the following results:

\smallskip
{\em
\begin{itemize}
\item[(i)] 
If  $\rho_0$ and $\rho_1$ are such that $u_1\le c_1$, then the
supersonic/subsonic regular reflection-diffraction solution
exists for each wedge angle $\theta_w\in (\theta_w^d, \frac{\pi}{2})$;

\smallskip
\item[(ii)]
If  $\rho_0$, $\rho_1$ are such that $u_1> c_1$, then there exists 
$\theta_w^a\in [\theta_w^d, \frac{\pi}2)$ such that
the regular reflection solution exists for
each  wedge angle $\theta_w\in (\theta_w^a, \frac{\pi}{2})$.
Moreover, if $\theta_w^a>\theta_w^d$, then, for the wedge angle $\theta_w=\theta_w^a$,
there exists an {\it attached} solution, {\it i.e.},
a solution of Problem {\rm 2} with  $\PtLwL=\PtLwR$.
\end{itemize}
The type of regular reflection (supersonic, Fig. {\rm \ref{fig:RegularReflection}}, or  subsonic,
Fig. {\rm \ref{fig:SubsoncRegularReflection}})
is determined by the type of state {\rm (2)} at $\PtIncW$.
The reflected shock $\PtIncW\PtLwL$ is $C^{2,\alpha}$-smooth.
The solution $\varphi$ is $C^{1,1}$ across the sonic arc for the supersonic
reflection-diffraction solution, which is optimal.
For the subsonic and sonic reflection-diffraction case
(Fig. {\rm \ref{fig:SubsoncRegularReflection}}), the solution is $C^{2,\alpha}$ in $\Omega$ near $P_0$.
Furthermore, the regular reflection-diffraction solution tends to the unique normal
reflection, when the wedge angle $\theta_w$ tends to $\frac{\pi}{2}$.
}

\medskip
These results are obtained by solving a free boundary problem for $\varphi$ in
$\Omega$.
The free boundary is the curved part of the reflected shock:
$P_1P_2$ on Fig. \ref{fig:RegularReflection}, and $P_0P_2$ on Fig. \ref{fig:SubsoncRegularReflection}.
We seek a solution $\varphi$ of equation (\ref{2-1}) in $\Omega$, satisfying the Rankine-Hugoniot conditions
(\ref{RH}) on the free boundary,
$\varphi_\nu=0$ on the boundary of the wedge and on $P_2P_3$,
and
$$
\varphi=\varphi_2, \quad D\varphi=D\varphi_2
$$
on $P_1P_4$ in the supersonic
case as shown in Fig. \ref{fig:RegularReflection}
and at $P_0$ in the subsonic case as shown in Fig. \ref{fig:SubsoncRegularReflection}.
The equation is expected to be elliptic in $\Omega$ for $\varphi$ (unknown {\it a priori} before
obtaining the solution), which is a part of the results.

\smallskip
To solve this free boundary problem, we define a class of admissible solutions of Problem 2,
which are solutions $\varphi$ with weak regular reflection-diffraction configuration,
such that, in the supersonic reflection case,
equation (\ref{2-1}) is strictly elliptic for $\varphi$
in $\overline\Omega\setminus \PtUpL\PtUpR$,
inequalities $\varphi_2\le\varphi\le \varphi_1$ hold in $\Omega$,
and the following monotonicity properties hold:
$$
\partial_\eta(\varphi_1-\varphi)\le 0, \quad D(\varphi_1-\varphi)\cdot e\le 0 \qquad \mbox{in $\Omega$},
$$
where $e=\frac{{P_0P_1}}{|P_0P_1|}$.
In the subsonic reflection
case, admissible solutions are defined similarly, with changes
corresponding to the structure of
subsonic reflection solution.
We derive uniform {\it a priori} estimates for admissible solutions with
wedge angles $\theta_w \in [\theta_w^d+\varepsilon, \frac\pi 2]$
for each $\varepsilon>0$,
and then apply the degree theory to obtain
the existence for each $\theta_w \in [\theta_w^d+\varepsilon, \frac\pi 2]$
in the class of admissible solutions,
starting from the unique solution for $\theta_w=\frac\pi 2$, called the normal reflection.
To derive the {\it a priori} bounds, we first obtain the estimates related to the geometry of the shock:
Show that the free boundary has a uniform positive distance from the sonic circle
of state (1) and from the wedge boundary away from $\PtLwL$ and $P_0$.
This allows to estimate the ellipticity of (\ref{2-1}) for $\varphi$ in $\Omega$
(depending on the distance to the sonic arc $P_1P_4$ for the supersonic reflection
and to $\PtIncW$ for the subsonic reflection).
Then we obtain the estimates near $P_1P_4$ (or $\PtIncW$
for the subsonic reflection) in scaled and weighted $C^{2,\alpha}$ for $\varphi$
and the free boundary, considering separately
four cases depending on $\frac{D\varphi_2}{c_2}$ at $\PtIncW$:

\medskip
\begin{itemize}
\item[(i)]
$\frac{|D\varphi_2(\PtIncW)|}{c_2} \ge 1+\delta$;

\smallskip
\item[(ii)]
 Supersonic (almost sonic): $1+\delta>
\frac{|D\varphi_2(\PtIncW)|}{c_2} > 1$;

\smallskip
\item[(iii)]
Subsonic (almost sonic): $1\ge 
\frac{ |D\varphi_2(\PtIncW)|}{c_2} \ge 1-\delta$;

\smallskip
\item[(iv)]
Subsonic: $
\frac{|D\varphi_2(\PtIncW)|}{c_2} \le 1-\delta$.
\end{itemize}

\medskip
In cases (i)--(ii), equation (\ref{2-1}) is degenerate elliptic in $\Omega$
near  $\PtUpL\PtUpR$ on Fig. \ref{fig:RegularReflection}.
In case (iii), the equation is uniformly elliptic in $\overline\Omega$,
but the ellipticity constant is small near $\PtIncW$ on Fig. \ref{fig:SubsoncRegularReflection}.
Thus, in cases (i)--(iii), we use the local elliptic degeneracy,
which allows to find a comparison function in each case, to show
the appropriately fast decay of $\varphi-\varphi_2$ near $\PtUpL\PtUpR$ in cases (i)--(ii) and near $\PtIncW$ in case (iii);
furthermore, combining with appropriate local non-isotropic rescaling to obtain the uniform ellipticity,
we obtain the {\it a priori} estimates in the weighted and scaled $C^{2,\alpha}$--norms, which are different in each
of cases (i)--(iii),
but imply
the standard $C^{1,1}$--estimates in cases (i)--(ii), and the standard
$C^{2,\alpha}$--estimates in case (iii).
This is an extension of methods of our earlier work \cite{Chen-F-6}.
In the uniformly elliptic case (iv), the solution is of subsonic reflection-diffraction
configuration as shown in Fig. \ref{fig:SubsoncRegularReflection},
and the estimates are more technically challenging than in cases (i)--(iii),
due to the lower {\it a priori} regularity
of the free boundary and since the uniform ellipticity does not allow
a comparison function that shows the decay of $\varphi-\varphi_2$ near $\PtIncW$.
Thus, we prove $C^{1,\alpha}$--estimates of $D(\varphi-\varphi_2)$ near $P_0$,
which imply the $C^{2,\alpha}$--estimates near $P_0$ for $\varphi$.
With all of these, we provide a solution to the von Neumann's conjectures.

\medskip
In Chen-Feldman-Xiang \cite{Chen-F-Xiang},
we have established the strict convexity of the curved (transonic)
part of shock
in the shock reflection-diffraction problem, as well as two other problems:
the Prandtl-Meyer reflection ({\it cf.} Bae-Chen-Feldman \cite{B-Chen-F-2013})
and shock diffraction ({\it cf.} Chen-Xiang \cite{chenxiang-shockdiffractionPotentialflow}).
In order to prove the convexity, we employ global properties of
admissible solutions,
including the existence of the cone of monotonicity discussed above.

\medskip
More details can be found in Chen-Feldman \cite{Chen-F-7}.

\vspace{.28in}
\noindent
{\bf Acknowledgements:} The research of
Gui-Qiang Chen was supported in part by
the UK EPSRC Science and Innovation
Award to the Oxford Centre for Nonlinear PDE (EP/E035027/1),
the NSFC under a joint project Grant 10728101, and
the Royal Society--Wolfson Research Merit Award (UK).
The work of Mikhail Feldman was
supported in part by the National Science Foundation under Grant
DMS-1101260
and by the Simons Foundation through
the Simons Fellows Program.

\bigskip

\end{document}